\documentclass[12pt]{amsart}

\usepackage{amscd,amsmath}
\usepackage{mathrsfs}
\usepackage{amsfonts}
\usepackage{amssymb}
\usepackage{enumerate}
\usepackage{setspace}
\usepackage{color}
\usepackage{nicefrac}
\usepackage{xfrac}
\usepackage{mathtools}

\usepackage[colorlinks=true,
            linkcolor=red,
            urlcolor=blue,
            citecolor=gray]{hyperref}

\usepackage{amsthm}
\usepackage{amssymb,latexsym,graphics,enumerate}
\usepackage[mathscr]{eucal}
\usepackage{amsmath,amsfonts,amsthm,amssymb}
\usepackage[left=1in,right=1in,top=1in,bottom=1in]{geometry}
\usepackage{tikz} 

\newcommand{\RR}{\mathbb{R}}
\newcommand{\TT}{\mathbb{T}}

\newcommand{\be}{\begin{equation}} 
\newcommand{\bea}{\begin{eqnarray}}
\newcommand{\ee}{\end{equation}} 
\newcommand{\eea}{\end{eqnarray}}
\newcommand{\ba}{\begin{aligned}}
\newcommand{\ea}{\end{aligned}}

\newcommand{\na}{\nabla}

\newcommand{\pa}{\partial}

\newtheorem{theorem}{Theorem}[section]

\newtheorem{lemma}{Lemma}[section]
\newtheorem{assumption}{Assumption}[section]

\newtheorem{remark}{Remark}[section]

\renewcommand{\geq}{\geqslant}
\renewcommand{\ge}{\geqslant}
\renewcommand{\leq}{\leqslant}
\renewcommand{\le}{\leqslant}
\def\d{{\textnormal{d}}}

\DeclareMathOperator{\trace}{\textnormal{trace}}

\def\u{{\mathbf u}}
\def\bu{{\mathbf u}}
\def\v{{\mathbf v}}
\def\x{{\mathbf x}}
\def\bx{{\mathbf x}}
\def\y{{\mathbf y}} 
\def\by{{\mathbf y}} 
\def\z{{\mathbf z}}

\def\vp{{\mathbf v}'}

\def\bw{{\mathbf w}}
\def\bff{{\mathbf f}}
\def\bF{{\mathbf F}}

\def\a{\phi}
\def\presure{{\mathbb P}}

\def\lam{\lambda}
\def\lamin{\lam_-} 
\def\lamax{\lam_+} 
\def\mumax{\mu_+} 

\def\hf{\tfrac{1}{2}}
\def\dv{\d\v}

\def\dx{\d\x}

\def\dy{\d\y}
\def\dz{\d\z}
\def\dt{{\d}t}

\def\etac{\eta_c}
\def\etar{c_*} 

\def\Cplus{C_+}
\def\Cminus{C_-}

\def\w{{\mathbf w}}
\def\z{{\mathbf z}}

\def\diam{D}

\def\supp{\textnormal{supp}}
\def\suppr{\supp\{\rho(t,\cdot)\}}
\def\rhoav{\phi*\rho} 
\def\rhoavpsi{\psi*\rho} 
\def\amp{\tau}
\def\ampone{\amp}
\def\ampzero{\amp m_0} 
\def\ddt{\frac{\textnormal{d}}{\textnormal{d}t}}
\def\nablaS{\nabla_{\!{}_S}}
\def\E{E} 
\def\kE{{\mathscr E}} 
\def\Etotal{{\mathscr E}} 
\def\delE{\delta\Etotal}

\def\amptwo{\amp}  

\numberwithin{equation}{section}

\begin{document}

\title[Entropy decrease and emergence of order in collective dynamics]{Entropy decrease and emergence of order\\in collective dynamics}

\author{Eitan Tadmor}
\address{Department of Mathematics and IPST, University of Maryland, College Park.}
\email{tadmor@umd.edu}

\date{\today}



\thanks{I thank Thomas Chen for helpful discussions. Research was supported  by  ONR grant N00014-2412659, NSF grant DMS-2508407 and  by the Fondation Sciences Math\'ematiques de Paris (FSMP) while being hosted by
LJLL at Sorbonne University.}
\date{\today}

\dedicatory{In memory of Ha\"{i}m Brezis}

\begin{abstract}
We study the hydrodynamic description of collective dynamics driven by velocity \emph{alignment}. 
  It is known that  such Euler alignment systems must flock towards  a limiting ``flocking'' velocity, provided their solutions remain globally smooth. 
   To address this question of global existence we proceed in two steps. (i) Entropy and closure. The system lacks a closure, reflecting lack of  detailed energy balance in collective dynamics. We discuss the decrease of entropy and the asymptotic behavior towards a mono-kinetic closure; and (ii) Mono-kinetic closure. We prove  that global regularity persists for all time for a large class of   initial conditions  satisfying a critical threshold  condition, which is intimately linked to the decrease of entropy.  The result applies in any number of spatial dimensions, thus addressing the  open question of existence beyond two dimensions.
\end{abstract}

\maketitle
\setcounter{tocdepth}{2}
\tableofcontents
\section{Introduction and statement of main results}
 We are concerned with the existence and large-time behavior of global smooth solutions  for the multi-dimensional Euler alignment system
\bea\label{eq:EA}
\quad \left\{
\begin{aligned}
&\rho_t+\na_\x\cdot(\rho \u)=0,\\
&(\rho \u)_t+\na_\x\cdot(\rho \u\otimes\u+\presure)=\amp\int \a(\x,\y)\big(\u(t,\y)-\u(t,\x)\big)\rho(t,\x)\rho(t,\y)\dy.
\end{aligned} \right.
\eea
A  solution pair of density-velocity, $(\rho,\u):\RR_+\times \Omega\rightarrow \RR_+\times \RR^n$,  is sought subject to compactly supported initial data 
$\big(\rho(0,\x),\u(0,\x)\big)=\big(\rho_0(\x),\u_0(\x)\big)$,  
either in the whole space, $\Omega=\RR^n$,  or in the $n$-dimensional torus $\Omega=\TT^n$.
System \eqref{eq:EA} is the large-crowd hydrodynamic description of `social agents' identified by positions and velocities, $\{\bx_i(t),\v_i(t)\}_{i=1}^N, \ N\gg1$,  governed by the Cucker-Smale alignment model \cite{CS2007}
\bea\label{eq:CS}
\left\{\begin{array}{rr}
\begin{aligned}
&\dot{\x}_i=\v_i\\
&\dot{\v}_i=\frac{\ampone}{N}\sum_{j=1}^N \phi(\x_i,\x_j)(\v_j-\v_i).
\end{aligned}\end{array}\right.
\eea
Alignment dynamics is a canonical model governing  emergence phenomena in collective dynamics  of flocking, swarming etc. The  dynamics is driven by a non-negative symmetric \emph{communication kernel}, $\phi(\x,\y)=\phi(\y,\x)\geq 0$, with amplitude $\amp>0$.  We have two main examples of symmetric kernels in mind  --- the canonical Cucker-Smale class of \emph{metric kernels}, see \cite{CS2007} 
\begin{equation}\label{eq:metricphi}
\phi(\x,\y)=\phi(|\x-\y|),
\end{equation}
 and the class of topologically-based kernels  \cite{ST2020b}
 \[
 \phi(\x,\y)=\phi_1(|\x-\y|)\phi_2(d_\rho(\x,\y)),\qquad d_\rho(\x,\y):=\int \limits_{{\mathcal C}(\x,\y)}\!\!\!\rho(t,\z)\dz,
 \] 
 reflecting the dependence on the mass   in an intermediate  `domain of communication', ${\mathcal C}(\x,\y)$, enclosed between $\x$ and $\y$.  
 
The Euler alignment system  \eqref{eq:EA} involves the pressure $\presure$ --- a symmetric positive definite tensor which should encode the thermodynamics of large-crowd collective dynamics. But \eqref{eq:EA}    is not a closed system: it   lacks closure of $\presure$ in terms of the macroscopic variables $\rho$ and $\u$.
 This reflects the fact that unlike physical particles,  the social agents engaged in  collective dynamics  form  a thermo-dynamically open system which is far from equilibrium and  does not admit a universal closure. 
 Indeed, most of the relevant literature \emph{assumes} a mono-kinetic closure, $\presure\equiv 0$. Accordingly, our study of solutions of \eqref{eq:EA} proceeds in two main stages:\newline
 (i) We investigate  a rather general class of so-called \emph{entropic pressures} introduced in \cite{Tad2023} and conclude that their strong solutions must approach mono-kinetic closure. Indeed, strong Euler alignment solutions with isentropic closure experience a uniform entropy \emph{decrease} to $- \infty$, and if we reject such a scenario then we must impose \emph{mono-kinetic closure},  $\presure\equiv 0$.\newline
  (ii) We consider  the pressure-less or mono-kinetic Euler alignment, proving existence of global smooth solutions under certain \emph{threshold conditions}. This  addresses the open question of existence for $n\geq 3$, extending the known results for dimensions $n=1,2$ \cite{TT2014,CCTT2016,HT2017}. 
 
 \subsection{The road to mono-kinetic closure}   
 To  investigate this mono-kinetic assumption we  let $\rho E$ denote the (total)   energy associated with the pressure in \eqref{eq:EA} (here and below, $|\w|$ denotes the $\ell_2$-norm of $\w$ and $|\w|_\infty$ denote its $L^\infty$-norm), 
 \be\label{eq:internal}
 \rho E:= \hf\rho|\u|^2 + \rho e,\qquad \rho e:=\hf \textnormal{trace}(\presure),
 \ee
Since the social agents engaged in collective dynamics  are often driven by energy  received from the ``outside'',  the detailed energy balance associated with \eqref{eq:EA} may be less relevant,  \cite[\S1.1]{VZ2012}.
Instead, lack of thermal equilibrium  in the form of  closure \emph{equalities}, can be relaxed  to certain \emph{inequalities}, which are compatible with the decay of \emph{energy fluctuations}. We impose the notion of an \emph{entropic pressure} in which   $\presure$, augmented with arbitrary ``heat-flux'' vector function ${\mathbf q}$, are required to satisfy the inequality, 
\begin{equation}\label{eq:entropy-ineq-E}
\begin{split}
(\rho \E)_t + \nabla_\x\cdot &(\rho \E\u + \presure\u +{\mathbf q}) \\
 & \leq -\amp\int\phi(\x,\y)\big(2\rho\E(t,\x)-\rho\u(t,\x)\cdot \u(t,\y)\big)\rho(t,\y)\dy.
 \end{split}
\end{equation}
This can be   expressed in an  equivalent form\footnote{The formal manipulation \eqref{eq:EA}${}_2\cdot \u -$\eqref{eq:EA}${}_1\times\frac{|\u|^2}{2}$ yields
\[
\big(\hf\rho|\u|^2\big)_t + \nabla_\x\cdot(\hf \rho|\u|^2\u+\presure\u)-\textnormal{trace}(\presure\nabla\u)=-\amp\int \phi(\x,\y)\rho\u(t,\x)\cdot\big(\u(t,\x)-\u(t,\y)\big)\rho(t,\y)\dy.
\]
Combining this equality with \eqref{eq:Int-entropic-e} is equivalent  with 
\eqref{eq:entropy-ineq-E}.}  in term of an entropic inequality for the internal energy, $\rho e$, stating  \cite[Definition 1.1]{Tad2023},
\begin{equation}\label{eq:Int-entropic-e}
(\rho e)_t + \nabla_\bx\cdot(\rho e \u+{\mathbf q}) + \textnormal{trace}\big(\presure\nabla\u\big) \leq -2\amp \int \phi(\x,\y)\rho e(t,\x)\rho(t,\y)\dy.
\end{equation}
 The entropic inequalities   \eqref{eq:entropy-ineq-E} or \eqref{eq:Int-entropic-e} are flexible enough to cover  pressure tensors derived from an underlying kinetic formulation discussed in \S\ref{sec:kinetic}, as well as the mono-kinetic closure $\presure\equiv 0$.
They imply \emph{depletion of   fluctuations}   and in particular, as we shall see later on, that the pressure in entropic alignment dynamics  approaches the mono-kinetic closure. In \S\ref{sec:mono-kinetic-decay} we discuss several results which quantify  this statement of ``approach towards mono-kinetic closure''. To state one  of these results, we need the following notations.
 
 \smallskip\noindent
 {\bf Notations}. We let $\nablaS\u$ denote the symmetric gradient, 
 $\displaystyle (\nablaS\u)_{ij}=\hf\Big(\frac{\partial u_i}{\partial x_j} + \frac{\partial u_j}{\partial x_i}\Big)$. We order the real   eigenvalues of $\nablaS\u$: \ $\lamin:=\lam_1\leq \lam_2 < \ldots \leq \lam_n=:\lamax$,  and in general, given a closed set $X\subset \RR $ we let $X_+$ and $X_-$ denote its largest, respectively smallest elements. We use $\etac$ to denote different positive constants. Finally, we use $\rhoav$ to denote the average density, or \emph{thickness}, 
 \[
  \rhoav(t,\x):=\int \phi(\x,\y)\rho(t,\y)\dy;
  \]
  of course, in case of metric kernel, this coincides with the usual notion of convolution. In \S\ref{sec:mono-kinetic-decay} we prove the following.
 \begin{theorem}[{\bf Decay towards mono-kinetic closure}]\label{thm:threshold}
 Let $(\rho, \u,\presure)$ be a strong entropic solution of \eqref{eq:EA} and assume  the following threshold condition holds: there exists a constant, $\etac>0$, such that 
 \begin{equation}\label{eq:Int-threshold}
 \lambda_-(\nablaS\u)(t,\cdot) + \amp\rhoav(t,\cdot)\geq\etac>0.
 \end{equation}
 Then $\displaystyle \int \|\presure(t,\x)\|\dx \leq e^{-\etac t}\int \|\presure_0(\x)\|\dx$. 
 \end{theorem}
 
\subsection{Communications kernels and thickness}
Different classes of communication kernels $\phi(\x,\y)$  are treated in the literature, classified according to their short- and long-range behavior where  $\x\approx \y$, and, respectively, $|\x-\y|\gg 1$.
The class of singular kernels $\phi(\x,\y)\sim |\x-\y|^{-\alpha}$ was treated in \cite{FP2024}
where it was shown that the enstrophy associated with strongly singular kernels, $\alpha\in[n,n+2)$, enforced mono-kinetic closure. Existence and flocking behavior of the  1D mono-kinetic case with strongly singular kernels was studied in \cite{DKRT2018, ST2017a,ST2017b,ST2018}.\newline 
In this paper we restrict attention to the case of bounded kernels, $\phi(\cdot,\cdot)\leq \phi_+$.
The results can be extended to the case of integrable kernel with weak singularity, $\phi(\x,\y)\sim |\x-\y|^{-\alpha}$ with $\alpha<n$. Alignment with such weakly singular kernels were studied in \cite{MMPZ2019}.

Next, we make  the distinction between long-range and short-range kernels.
Assume that $\phi(\cdot,\cdot)$ admits a metric lower envelope,
$\phi(\x,\y)\gtrsim \varphi(|\x-\y|)$, with decreasing radial $\varphi$.  
Heavy-tailed kernels is the subclass of long-range kernels that naturally arise in  connection with the unconditional flocking behavior of \eqref{eq:EA}, \cite{HT2008,HL2009,CFTV2010,CFRT2010}
\[
\int^\infty \varphi(r) \d r =\infty.
\]
Heavy-tails rule out short-range kernels  with  finite support which are important in applications. In this latter case, flocking  is  secured if the Euler alignment dynamics remains uniformly non-vacuous, $\rho\geq \rho_->0$ \cite[Theorem 3]{Tad2021}, \cite[\S7.3]{Shv2024}.

\noindent
We proceed by making  the following thickness condition. This corresponds to the notion of ball-thickness introduced in \cite[\S3.7.2]{Shv2024}.
The scenario of (uniform) thickness covers both cases of long-range (-- heavy-tailed) kernels and short-range (non-vacuous) kernels; this will be discussed in \S\ref{sec:mono-kinetic-existence} below.

\begin{assumption}[{\bf  Thickness}]
The following thickness condition  holds
\[
 \int^\infty \min_\x \rhoav(t,\x)\dt =\infty.
\]
In particular, a \emph{uniform thickness} condition holds 
if there exists a constant, $\etar>0$, such that
\be\label{eq:Int-thickness}
\rhoav(t,\cdot) \geq \etar>0.
\ee
\end{assumption}

\subsection{Decrease of  entropy and emergence of order} We restrict attention to scalar pressure,  $\presure=p{\mathbb I}_{n\times n}$. In this case we have the isentropic closure, $p= \frac{1}{n} \textnormal{trace}(\presure) =  \frac{2}{n}\rho e$,
and the entropic inequality \eqref{eq:Int-entropic-e} reads
\[
p_t + \u\cdot \nabla_\x p +\gamma p \nabla_\bx\cdot\u \leq -2\amp p \rhoav, \qquad \gamma :=1+\frac{2}{n}.
\] 
Next we manipulate --- multiplying the last inequality by $\rho^{-\gamma}$ and adding a multiple of the mass equation,  $-\gamma \rho^{\gamma-1}p \times \eqref{eq:EA}_{1}$, to find
\[
(p\rho^{-\gamma})_t + \u\cdot\nabla_\x(p\rho^{-\gamma}) \leq -2\amp p\rho^{-\gamma}\rhoav.
\]
We conclude that the entropy inequality we imposed in \eqref{eq:Int-entropic-e} amounts to an inequality on the  specific entropy, $\ln(p\rho^{-\gamma})$,
\begin{equation}\label{eq:Int-S-decreases}
S_t +\u\cdot\nabla_\x S\leq -2\amp\rhoav < 0, \qquad S:=\ln(p\rho^{-\gamma}).
\end{equation}
The conclusion that entropy  decreases in time is  quite striking in the sense  
that it  counters  the usual tendency of entropy to quantify disorganization as expressed by the second law of   entropy increase  in physical systems. 
Ever since  \cite{Sch1944} stipulated   that   life maintains its own organization by extracting order from its environment,  stating  ``\emph{life ... feeds on negative entropy}'', there have been many attempts to link living systems to a reversed second law which drives organisms towards higher ordered content, \cite{Cha1978, Ave2003}. 
We have shown here that the self-organization of alignment dynamics decreases the entropy by communicating ``information'' in the whole
flow-field  through decrease of fluctuations. In \S\ref{sec:mono-kinetic-decay} we  further elaborate on implications of the entropy decrease in alignment dynamics. 
Here is one manifestation of  \eqref{eq:Int-S-decreases} proved in \S\ref{sec:entropy}.
\begin{theorem}\label{thm:dichotomy}
Let $(\rho, \u)$ be a strong entropic solution of \eqref{eq:EA} with isentropic closure\newline $\presure= \frac{2}{n}\rho e{\mathbb I}_{n\times n}$, and assume the   thickness condition holds. Then  there is an  entropy decay 
$S(t,\cdot) \stackrel{t\rightarrow\infty}{\longrightarrow} -\infty$,
 unless $\presure\equiv 0$.
\end{theorem}
Theorem \ref{thm:dichotomy} expresses the following dichotomy between two possible scenarios for a non-vacuous strong
 isentropic solutions: either, for $\rho, p>0$,  there is large time decay towards a mono-kinetic closure --- in fact the uniform thickness \eqref{eq:Int-thickness} implies exponential decay  \emph{uniformly for all} $\x$
 \begin{equation}\label{eq:exp-decay-ent}
 \max_\x p\rho^{-\gamma}(t,\x) \lesssim e^{-2\amp \etar t} \stackrel{t\rightarrow\infty}{\longrightarrow}0,
 \end{equation}
or else a mono-kinetic closure $p\equiv 0$. If we reject the first scenario then we conclude Euler alignment system can admit global strong solutions only in the setting of mono-kinetic closure.
 \subsection{Strong solutions with mono-kinetic closure} We now turn  to investigate the existence of strong solutions of the alignment dynamics system \eqref{eq:EA} with mono-kinetic closure $\presure\equiv 0$,
\bea\label{eq:Int-mono-kinetic}
\quad \left\{
\begin{aligned}
&\rho_t+\na_\x\cdot(\rho \u)=0,\\
&\u_t+\u\cdot\na_\x\u=\amp\int \a(\x,\y)\big(\u(t,\y)-\u(t,\x)\big)\rho(t,\y)\dy.
\end{aligned} \right.
\eea

The 1D system \eqref{eq:Int-mono-kinetic} with metric kernel \eqref{eq:metricphi} admits global smooth solution if and only if the initial condition satisfies the following lower critical threshold  \cite{CCTT2016}, $u'_0+\amp\rhoav_0>0$.
This was extended to the class of uni-directional flows,  $\u(t,\x) := (u(t,\x), 0,\ldots,0)$ with $u: \RR_+\times \Omega \mapsto \RR$, in
 \cite{LS2020}, and to  2D metric case in \cite{HT2017}. In \S\ref{sec:mono-kinetic-existence} we address the open question of existence of solution of Euler alignment system \eqref{eq:Int-mono-kinetic} with general symmetric kernels in $n>2$ spatial dimensions. The case of metric kernels is summarized in the following.
\begin{theorem}[{\bf Global strong solutions with sub-critical data}]\label{thm:CS-metric}
\noindent
Consider the mono-kinetic Euler alignment  \eqref{eq:Int-mono-kinetic}
with  metric communication kernel $\phi(\x,\y)=\phi(|\x-\y|)>0$ satisfying 
a uniform thickness \eqref{eq:Int-thickness}, $\rhoav(t,\cdot)\geq \etar >0$. 
The system is  subject to initial conditions $(\rho_0, \u_0) \in  L_+^1\times W^{1,\infty}$ with velocity fluctuations which do not exceed,
\begin{equation}\label{eq:Int-V0}
8 |\phi'|_\infty \cdot \delta u_0 < \amp\etar^2, \qquad \delta u(t):=\sup_{{\x,\y\in \suppr}}\hspace*{-0.2cm} |\u(t,\x)-\u(t,\y)|. 
\end{equation}
Assume that the initial data $(\rho_0,\u_0)$ satisfy a sub-critical threshold condition 
\begin{equation}\label{eq:Int-CTnD}
\lamin(\nablaS\u_0)(\x) +\amptwo\rhoav_0(\x) \geq \etac>0, \quad \etac:=\tfrac{1}{2}\amp\etar.
\end{equation}
Then the Euler alignment  system \eqref{eq:Int-mono-kinetic}--\eqref{eq:Int-CTnD}  
 admits a global smooth solution, $(\rho(\cdot,t),\u(\cdot,t))\in L^1_+\times W^{1,\infty}(\RR^n)$, with uniformly bounded velocity gradient, 
$|\nabla \u(t,\cdot)|_{L^\infty} \leq \max\left\{|\nabla \u_0|_{L^\infty},\etar,C_0\right\}<\infty$.
 \end{theorem}
 
 \begin{remark}[{\bf What does the threshold condition mean?}]\label{rmk:what}
 Considering the limiting case $\amp=0$ then the threshold condition \eqref{eq:Int-CTnD} requires $\lamin(\nablaS\u_0)(\x)>0$. In the one-dimensional case, this reflects the fact that the inviscid Burgers' equation  admits global smooth solution for increasing profile $u'_0>0$; however, this is  a rather restricted set of initial profiles. Similarly, in the $n$-dimensional case, the threshold  \eqref{eq:Int-CTnD} with $\amp=0$ reflects  global smooth solutions of the pressure-less Euler $\u_t+\u\cdot\nabla_\x\u=0$ for the restricted class of initial configurations $\lamin(\nablaS\u_0)(\x)>0$ (for which $\nabla_\x\cdot\u_0>0$ excludes, for example, $|\u(t,\cdot)|\stackrel{|\x|\rightarrow\infty}{\longrightarrow}0$). Thus, the essence of the threshold condition \eqref{eq:Int-CTnD} is securing global existence for a large set of initial configurations, by  allowing $\lamin(\nablaS\u_0)$ (or $\nabla_\x\cdot \u_0$) to admit negative values dictated by the local thickness, $\amp\rhoav$. Observe that  according \eqref{eq:Int-CTnD}, the admissible  values of $\lamin(\nablaS\u_0)$ include the negative range  
 \[
 \lamin(\nablaS\u_0) \geq \etac - \amp \rhoav, \qquad  \etac - \amp \rhoav\leq -\hf\amp\etar, \ \  \etac=\hf\amp\etar.
 \]
  This type of critical threshold which secures global regularity with negative  ''initial slopes'' as long as they are  ``not too negative'' was introduced in the context of  Euler-Poisson equations in \cite{ELT2001,LT2002} and is found useful for Euler alignment models, \cite{TT2014,CCTT2016}; see \cite{Shv2024} and the references therein.
 \end{remark}
  \begin{remark}[{\bf Comparing the threshold conditions}] It is instructive that the  threshold condition for global existence sought in \eqref{eq:Int-CTnD}, $\lamin(\nablaS\u_0) +\amptwo\rhoav_0 \geq \etac>0$, is the same condition we met earlier, \eqref{eq:Int-threshold}, in the context of decay towards mono-kinetic closure.\newline
  We compare this threshold condition with the results available in current literature on  the global regularity of mono-kinetic Euler alignment system \eqref{eq:Int-mono-kinetic} in one- and two dimensions. 
Global smooth solutions in the 1D case, and in the more general setup of \emph{uni-directional} flows,  exist  if and only if the initial configuration satisfies the initial threshold $u'_0(x)+\amp\rhoav_0(x)\geq 0$, \cite{CCTT2016,LS2020,Les2020}. This corresponds to the limiting case $\etac=0$.
For the role of the zero set $\{x \ | \ u'_0(x)+\amp\rhoav_0(x)= 0\}$ we refer to \cite{LLST2022}.\newline
 A sufficient  threshold for 2D regularity, \cite[Theorem 2.1]{HT2017},
requires the initial threshold $\nabla\cdot \u_0 + \amp\rhoav_0 >0$ \emph{and}
$(\lambda_+-\lambda_-)(\nablaS \u_0) \leq \amp\delta_0$ with $\delta_0=\hf m_0\phi(D_\infty)$. Noting that $\nabla\cdot\u=\lamax+\lamin$ then 
\[
\lamin(\nablaS\u_0)+\amptwo\rhoav_0  =\frac{\lamax+\lamin}{2}+\amptwo\rhoav_0-\frac{\lamax-\lamin}{2}\geq \frac{\amp\rhoav}{2} -\frac{\amp\delta_0}{2}, \quad \rhoav\geq m_0\phi(D_\infty).
\]
Thus, this 2D result is covered by  Theorem \ref{thm:CS-metric}  subject to  threshold $\lamin(\nablaS\u_0) +\amptwo\rhoav_0 \geq \etac>0$ with $\etac=\hf \delta_0$.
\end{remark}

\section{Kinetic formulation}\label{sec:kinetic}
The passage from the agent-based description \eqref{eq:CS} to the hydrodynamic description  \eqref{eq:EA} goes through a kinetic formulation,  
\begin{equation}\label{eq:Q}
\begin{split}
\partial_t f +\v\cdot\nabla_\x f = -\amp \nabla_\v\cdot Q_\phi(f,f), \qquad (t,\x,\v)\in \RR_+\times\Omega\times \RR^n,
 \end{split}
\end{equation}
which governs the empirical distribution
\[
f=f_N(t,\x,\v):=\frac{1}{N}\sum_{j=1}^N\delta_{\x-\x_j(t)}\otimes \delta_{\v-\v_j(t)},
\]
and is driven by  pairwise  communication protocol on the right of \eqref{eq:CS}${}_{2}$
\begin{equation}\label{eq:Qphi}
 Q_\phi(f,f):=\iint_{\RR^{n}\times \Omega}\phi(\x,\y)(\vp-\v)f(t,\x,\v)f(t,\y,\v')\dv'\dy.
\end{equation}
The formal derivation of \eqref{eq:Q} was introduced in \cite{HT2008} and was justified in increasing order of rigor in \cite{HL2009,CFTV2010,FK2019, NP2022,NS2022,PT2022}.
For large crowds of $N$ agents, $N\gg1$,  the dynamics \eqref{eq:Q}  is captured by its  first three  limiting moments which are assumed to exist in a proper sense: the density $\displaystyle \rho(t,\x):= \lim_{N\rightarrow \infty}\int  f_N(t,\x,\v)\dv$,  momentum $\displaystyle \rho \u(t,\x) :=\lim_{N\rightarrow \infty} \int \v f_N(t,\x,\v)\dv$ 
and  pressure  $\displaystyle \presure(t,\x):=\lim_{N\rightarrow \infty}\int (\v-\u)\otimes(\v-\u) f_N(t,\x,\v)\dv$, governed by  the Euler alignment system \eqref{eq:EA}.

\smallskip
To address the lack of closure we consider the  energy balance associated with  
  the limiting quadratic moment of \eqref{eq:Q} (which is assumed to exist)
$\displaystyle   \rho \E(t,\x):=\lim_{N\rightarrow \infty} \int \frac{1}{2}|\v|^2 f_N(t,\x,\v)\dv$, 
\[
(\rho \E)_t + \nabla_\x\cdot (\rho \E\u + \presure\u +{\mathbf q}) =-\amp\int\phi(\x,\y)\big(2\rho\E(t,\x)-\rho\u(t,\x)\cdot \rho\u(t,\y)\big)\dy.
\]
The total energy $\rho E$ admits a decomposition  into  its kinetic and internal parts,
$\rho \E = \hf \rho |\u|^2 + \rho e$,  corresponding to  the two-term decomposition\footnote{Here  and below, ``=''  is interpreted as ``equality modulo linear moments'', noting that linear moments vanish $\displaystyle \int (\v-\u)f_N(t,\x,\v)\dv=0$.} $|\v|^2$ {``=''} $|\u|^2+|\v-\u|^2$, namely, 
\[
\rho \E = \hf \rho |\u|^2 + \rho e, \qquad 
(\rho e)(t,\x):=\lim_{N\rightarrow \infty}\frac{1}{2}\int|\v-\u|^2 f_N(t,\x,\v)\dv=\hf \textnormal{trace}(\presure).
\]
The energy flux on the left, $\rho \E\u + \presure\u +{\mathbf q}$, is recovered as the quadratic  moment of \eqref{eq:Q} corresponding to the three-term decomposition, 
\[
|\v|^2\v \ \textnormal{``=''} \ |\v|^2\u +2(\v-\u)\otimes(\v-\u)\u + |\v-\u|^2(\v-\u).
\]
 This results in an energy flux expressed in terms of  the pressure   $\presure$ and ``heat flux''  $\displaystyle {\mathbf q}(t,\x):= \lim \limits_{N\rightarrow \infty}\frac{1}{2}\int (\v-\u)|\v-\u|^2f_N(t,\x,\v)\dv$.
To address  the lack of closure problem  we impose the notion of an entropic pressure which relaxes the energy equality, requiring    $(\presure,{\mathbf q})$ to satisfy the corresponding inequality \eqref{eq:entropy-ineq-E}
\[
(\rho \E)_t + \nabla_\x\cdot (\rho \E\u + \presure\u +{\mathbf q}) \leq -\amp\int\phi(\x,\y)\big(2\rho\E(t,\x)-\rho\u(t,\x)\cdot \u(t,\y)\big)\rho(t,\y)\dy.
\]

\section{Decay of fluctuations towards mono-kinetic closure}\label{sec:mono-kinetic-decay}
We discuss various scenarios in which an entropic pressure decays towards mono-kinetic closure.
\paragraph{{\bf Mono-kinetic closure with heavy-tailed kernels}} The formulation of the entropy inequality   in terms of the total energy is equivalent to imposing an instantaneous  decay of energy fluctuations. 
Indeed, integration of \eqref{eq:entropy-ineq-E} implies
\[
\begin{split}
\ddt \int \big(\hf &\rho |\u|^2(t,\x)  + \rho e(t,\x)\big) \dx \\
 & \leq -\amp \iint \phi(\x,\y)\big(\rho|\u|^2(t,\x)+2\rho e(t,\x)
-\rho\u(t,\x)\cdot\u(t,\y)\big)\rho(t,\y)\dy.
\end{split}
\]
Symmetrization of the integrand on the left using the fact  that $\displaystyle  \int \rho\u(t,\x) \dx \equiv m_0$, and symmetrization of the integrand of the right using the assumed symmetry of $\phi$, finally yields 
 the decay of fluctuations 
\begin{equation}\label{eq:Int-Efluc}
\begin{split}
 \ddt &\frac{1}{2m_0}\iint \Big(\hf|\u(t,\x)-\u(t,\y)|^2 +e(t,\x)+e(t,\y)\Big)\rho(\x)\rho(\y)\dx\dy \\
  & \leq -\amp \iint \phi(\x,\y)\Big(\hf|\u(t,\x)-\u(t,\y)|^2 + e(t,\x) +e(t,\y)\Big) \rho(\x)\rho(\y)\dx\dy.
  \end{split}
 \end{equation}
 We conclude the following result of \cite[Theorem 4.1]{Tad2023}.
 \begin{theorem}
 Consider the Euler alignment system with heavy-tailed kernel 
  \be\label{eq:heavy-tailed}
 \phi(\x,\y) \gtrsim \langle |\x-\y|\rangle^{-\theta}, \quad \theta\leq 1, \qquad \langle r\rangle:=(1+r^2)^{\nicefrac{1}{2}}
 \ee
  Let $(\rho,\u)$ be   a non-vacuous strong solution of \eqref{eq:EA}
 with entropic pressure $\presure$. Let $D(t)$ denote the diameter of $\textnormal{supp}\, \rho(t,\cdot)$, and assume the dispersion bound,
 \be\label{eq:dispersion}
 \int \langle D(t)\rangle^{-\theta}\dt =\infty, \qquad D(t):=\max_{\x,\y\in \textnormal{supp}\, \rho(t,\cdot)}|\x-\y|.
 \ee
 Then unconditional flocking holds 
 \[
 \iint \Big(\hf|\u(t,\x)-\u(t,\y)|^2 +e(t,\x)+e(t,\y)\Big)\rho(t,\x)\rho(t,\y)\dx\dy \stackrel{t\rightarrow \infty}{\longrightarrow} 0,
 \]
 and in particular,  asymptotic mono-kinetic closure holds, $\displaystyle \int \|\presure(t,\x)\|\dx \stackrel{t\rightarrow \infty}{\longrightarrow} 0$.
 \end{theorem}
This shows that under the dispersion bound \eqref{eq:dispersion}, the same mechanism that is responsible for unconditional flocking, $\|\u(t,\x)-\u(t,\y)\|_{L^2(\d\rho(\x)\d\rho(\y)} \stackrel{t\rightarrow \infty}{\longrightarrow} 0$, drives the alignment dynamics with \emph{any} entropic   pressure towards mono-kinetic closure, $\displaystyle \|\presure(t,\cdot)\|_{L^1}=\int \rho e(,\x) \dx   
\stackrel{t\rightarrow \infty}{\longrightarrow} 0$. 

\smallskip
\paragraph{{\bf Mono-kinetic closure with a threshold condition}}
The dispersion bound sought in \eqref{eq:dispersion} can be replaced by a threshold condition. To this we revisit the entropy inequality \eqref{eq:Int-entropic-e} which is re-written in  terms of the 
re-scaled pressure
$\displaystyle  \overline{\presure}= \frac{\presure}{2\rho e}$
\be\label{eq:Int-entropic-er}
\partial_t (\rho e) + \nabla_\x\cdot (\rho e\u+{\mathbf q}) \leq 
-\big(\trace(\overline{\presure}\nabla\bu)+\amp\rhoav\big) 2\rho e(t,\x), \quad 
    \overline{\presure}:=\frac{\presure}{2\rho e}
\ee
 Observe that 
 $\overline{\presure}$ is a symmetric positive definite matrix with unit trace. For such unit-trace matrices  we have 
\[
\trace(\overline{\presure}M) \geq \lamin(M_S), \qquad M_S:=\frac{1}{2}(M+M^\top).
\]
Indeed, if we  let $\{\lambda_i>0, {\mathbf w}_i\}$ be the complete eigen-system of $\overline{\presure}$ then
\[
\textnormal{trace}(\overline{\presure}M) =\sum_i \langle \overline{\presure}M{\mathbf w}_i,{\mathbf w}_i\rangle=
\sum_i \lambda_i(\overline{\presure})\langle M{\mathbf w}_i,{\mathbf w}_i\rangle
 \geq \sum_i \lambda_i(\overline{\presure})\lambda_-(M_S)
=\lambda_-(M_S).
\]
In particular, therefore, \eqref{eq:Int-entropic-er} yields
\[
\partial_t (\rho e) + \nabla_\x\cdot (\rho e\u+{\mathbf q}) \leq 
-2\big(\lamin(\nablaS  \u)+\amp\rhoav\big)\rho e(t,\x).
\]
We conclude that if  the threshold condition \eqref{eq:Int-threshold} holds 
\begin{equation}\label{eq:CT} 
\eta(t,\x)\equiv\eta(\rho(t,\x),\u(t,\x)):= \lamin(\nablaS \u)(t,\x) + \amp \rhoav(t,\x) \geq \etac >0,
\end{equation}
then the decay of  internal energy follows,
\[
\ddt \int \rho e (t,\x)\dx\leq -2\eta_c \int \rho e(t,\x)\dx,
\] 
which in turn proves Theorem \ref{thm:threshold}.
  The key question is whether  such threshold inequality $\eta(t,\x)\geq\eta_c>0$ persists in time.
  Observe that \eqref{eq:CT} is independent of the thermo-dynamical state $\{e, \presure\}$, and in particular, therefore, applies to the mono-kinetic closure $\presure=0$. This motivates our 
  search in \S\ref{sec:mono-kinetic-existence} for a \emph{critical threshold} in the space of initial configurations,   $ \eta(\rho_0(\x),\u_0(\x))\geq \eta_c>0$.

\smallskip
\subsection{Entropy decrease and mono-kinetic closure}\label{sec:entropy}
We restrict attention to the case of isentropic closure $p=(\gamma-1)\rho e$ which led to the reverse entropy inequality \eqref{eq:Int-S-decreases}.
\be\label{eq:S-decreases}
S_t +\u\cdot\nabla_\x S \leq -2\amp\rhoav<0, \qquad S=p\rho^{-\gamma}.
\ee
The Euler alignment system \eqref{eq:EA},\eqref{eq:entropy-ineq-E} can be viewed as hyperbolic system
of conservation laws for $\bw=(\rho,\rho \u,\rho E)$, which we abbreviate as 
\[
\bw_t +\textnormal{div}\,\bff(\bw) = \amp {\mathbb A}(\bw), \quad \amp>0
\] 
where $\bff(\bw)$ is the flux and ${\mathbb A}(\bw)$,
\[
{\mathbb A}(\bw):={\small \begin{bmatrix*}[l]\ \ \ \ 0 \\ \displaystyle  \ \ \ \int \phi(\x,\y)\big(\u(t,y)-\u(t,\x)\big)\rho(t,\x)\rho(t,\y)\dy \\ \displaystyle \leq \int \phi(\x,\y)\big(2\rho E(t,\x)\rho(t,\y)-\rho \u(t,\x)\cdot\rho\u(t,\y)\big)\dy\end{bmatrix*}},
\]
 encodes the alignment terms on the right of \eqref{eq:EA},\eqref{eq:entropy-ineq-E}. In this context, $(-\rho S, -\rho\u S)$ forms  an \emph{entropy pair}:   combining  \eqref{eq:S-decreases} and the mass equation, \eqref{eq:S-decreases}$\times \rho$ $+$ 
\eqref{eq:EA}$\mbox{}_1 \times S$,  implies  that 
 the convex entropy $-\rho S$ is \emph{increasing},
\[
U(\bw)_t +\nabla_\x\cdot\bF(\bw) >0, \qquad U(\bw)=-\rho S, \ \ \bF(\bw)=-\rho\u S,
\]
in contrast to the celebrated statement  of convex entropy decrease in presence of  diffusion, ${\mathbb D}(\bw)$, \cite{Lax1957,God1962},\cite[\S7]{Kru1970},\cite{Lax1971,Daf2005},
\[
\bw_t +\textnormal{div}\,\bff(\bw) = \sigma {\mathbb D}(\bw) \ \ 
\leadsto \ \  U(\bw)_t +\nabla_\x\cdot\bF(\bw) <0, \qquad \sigma>0.
\]
It follows  that alignment and diffusion compete in driving the dynamics in different directions of  increasing order and, respectively, disorder. This is realized in the reversed   entropy inequality \eqref{eq:S-decreases} which implies the maximum principle, $S(t,\cdot)\leq \max S_0$, in contrast to the
minimum principle in the vanishing difussive Euler equations \cite{Tad1986}
\[
S_t +\u\cdot\nabla_\x S\geq 0 \ \ \leadsto  \ \ S(t,\cdot)\geq \min S_0.
\] 
In fact, uniform thickness in \eqref{eq:S-decreases} implies 
\[
S^\eta_t +\u\cdot\nabla_\x S^\eta <0, \qquad S^\eta(t,\x):=S(t,\x)+2\amp\etar t,
\]
which in turn enforces the decay of $S(t,\cdot)$ asserted in Theorem \ref{thm:dichotomy}. We verify this by setting\footnote{$X^+:=\left\{\begin{array}{ll} X & X>0\\ 0 & X\leq 0\end{array}\right.$ denotes the ``positive part of $X$''.} ${\mathcal H}(S):= (S-S_*)^+$ with $S_*$ to be determined. Since ${\mathcal H}(\cdot)$ is non-decreasing then  $(\partial_t +\u\cdot\nabla_\x){\mathcal H}(S^\eta)\leq 0$, and  the entropy inequality follows
\[
(\rho {\mathcal H}(S^\eta))_t +\nabla_\x\cdot \big(\rho \u {\mathcal H}(S^\eta)\big) \leq 0.
\]
Now let $\displaystyle S_*:=\max_\x (S^\eta)_0= \max_\x S_0$, then integration of the entropy inequality yields 
\[
\int_{\x} \rho (S^\eta(t,\x)-S_*)^+\dx \leq \int_{\x} \rho (S_0(\x)-S_*)^+ \dx =0,
\]
implying that the non-negative integrand on the  left must vanish. Hence $\displaystyle S^\eta(t,\x)\leq S_* \leq \max_\x S_0$ and  \eqref{eq:exp-decay-ent}  follows,
\[
p\rho^{-\gamma}(t,\x) \leq K e^{-2\amp\etar t}, \quad K=e^{ \max_\x \!\!S_0}.
\]
Similar arguments apply  to the uniform decay $S(t,\cdot) \stackrel{t\rightarrow \infty}{\longrightarrow}-\infty$ asserted  in Theorem \ref{thm:dichotomy} under a general  thickness condition.

The ``competition'' between entropy production and entropy dissipation  is demonstrated in the one-dimensional  Navier-Stokes alignment (NSA) equations which we abbreviate
\begin{equation}\label{eq:NSA}
\bw_t +\bff(\bw)_x = \amp {\mathbb A}(\bw)+\sigma{\mathbb D}(\bw), \ \ \sigma{\mathbb D}(\bw):=
\sigma_1\begin{bmatrix*}[l] 0\\ u\\ \hf u^2\end{bmatrix*}_{xx} \!\!\!\!\!+ \sigma_2\begin{bmatrix*}[l] 0\\ 0\\ T\end{bmatrix*}_{x} \ \  C_v T = e,
\end{equation}
with  the corresponding entropy balance  
\begin{equation}\label{eq:NSA-entropy}
(-\rho S)_t+\big(-\rho u S +\sigma_2 (\ln T)_x\big)_x = 2\amp( \rhoav) \rho-\sigma_1\frac{u_x^2}{T}-\sigma_2 \left(\frac{T_x}{T}\right)^2.
\end{equation}
The entropy increase rate on the first term on the right of \eqref{eq:NSA-entropy} follows from \eqref{eq:S-decreases}; the terms involved the temperature $T$ encode the usual entropy decrease
associated with the diffusion term in Navier-Stokes, $\sigma{\mathbb D}(\bw)$. Observe that whenever the NSA dynamics \eqref{eq:NSA} with vanishing amplitudes $\amp,\sigma \ll1$  develops sharp gradients,  then the diffusive entropy dissipation terms,  $\sigma_1 u_x^2, \sigma_2 T_x^2 \gg 1$, dominate the bounded entropy increase due to alignment, $\rhoav \leq Const.$.

Finally, we note that the ``competition''  between entropy production and entropy dissipation  is  realized already at the kinetic level. We consider the kinetic formulation \eqref{eq:Q} driven by both alignment and diffusion,
 \cite[\S6]{Shv2024},\cite{Shv2025}
\[
\partial_t f +\v\cdot\nabla_\x f +\amp \nabla_\v\cdot Q_\phi(f,f)= \sigma\Delta_\v f.
\]
It follows that $H(f):=f\log f -f$ satisfies
\be\label{eq:kinetic-H}
\partial_t \int H(f)\dv +\nabla_\x \cdot \int \v H(f)\dv =\amp\phi*\rho \int f\dv -4\sigma  \int |\nabla_\v \sqrt{f}|^2 \dv.
\ee
The first term on the right shows that alignment increases the kinetic entropy  functional $\displaystyle \int H(f)\dv$  due to communication satisfying the   uniform thickness assumption \eqref{eq:Int-thickness},  $\amp\phi*\rho  f>\amp\etar f>0$. This reversed $H$ theorem was already observed by us in \cite[\S6]{HT2008}. In contrast, the  decrease of the kinetic entropy  functional due to diffusion is dictated by Fisher information $-4\sigma |\nabla_\v \sqrt{f}|^2<0$. They balance each other to zero entropy production with a Maxwellian profile
\[
f(t,\x,\v)=\frac{\rho}{(2\pi\theta)^{n/2}}e^{-\tfrac{|\v-\u|^2}{2\theta}}, \qquad
 \theta(t,\x) \sim \left(\frac{\sigma}{\amp}\frac{\rho}{\rhoav}\right)^{\nicefrac{(\gamma-1)}{\gamma}}.
\]

\section{Alignment with mono-kinetic closure}\label{sec:mono-kinetic-existence}
 Our  main result settles the open question of existence of strong solutions of the multi-dimensional mono-kinetic  Euler alignment system \eqref{eq:EA} in  $n\geq 3$ dimensions subject to  general $C^1$ symmetric communication kernels. 
\subsection{Existence of strong solutions}\label{sec:existence}
 The mono-kinetic closure reduces the momentum equation \eqref{eq:EA}$\mbox{}_2$ to \eqref{eq:Int-mono-kinetic}$\mbox{}_2$
 \be\label{eq:mono-momentum}
 \u_t(t,\x) +\u\cdot\nabla_\x\u(t,\x) \!=\! \int \!\!\phi(\x,\y)(\u(t,\y)-\u(t,\x))\rho(t,\y)\dy,\  \x\in \textnormal{supp}\,(\rho(t,\cdot).
 \ee

The existence result involves  three quantities   which  do not increase in time: the mass $\displaystyle m(t):=\int \rho(t,\x)\dx=m_0$, the  velocity fluctuation, $\delta u(t)$, (see \eqref{eq:max-principle} below),
\begin{equation}\label{eq:noincrease}
\delta u(t)  \leq \delta u_0,  \  \quad \ \ \delta u(t):=\sup_{\x,\y} \big\{|\u(t,\x)-\u(t,\y)|\, : \, {\x,\y\in \suppr}\big\},
\end{equation}
and the total kinetic energy, $\displaystyle \kE(t):= \int \hf \rho |\u|^2(t,\x)\dx\leq \kE(0)$.
We assume that their initial amplitudes  are not too large relative to the uniform thickness, $\etar>0$,
\begin{equation}\label{eq:V0}
(8\alpha_0+4\beta_0)m_0 < \amp\etar^2 , \qquad \left\{\begin{array}{ll}
\alpha_0:= |\nabla_\x\phi|_\infty |\delta u_0|_\infty\cr
 \beta_0:= \displaystyle |(\nabla_\x+\nabla_{\y})\phi|_\infty \left(\dfrac{2\kE_0}{m_0}\right)^{1/2}. \end{array}\right.
\end{equation}

\begin{theorem}[{\bf Global strong solutions with sub-critical data}]\label{thm:CS}
\noindent
Consider the multiD Euler alignment system \eqref{eq:Int-mono-kinetic} with $C^1$ symmetric communication kernel,  $\phi(\cdot,\cdot)>0$, satisfying the uniform thickness  \eqref{eq:Int-thickness}, $ \rhoav(t,\cdot)\geq \etar>0$, and subject to initial conditions $(\rho_0, \u_0) \in  L_+^1\times W^{1,\infty}$
with bounded amplitudes \eqref{eq:V0}.\newline
 Assume that the threshold condition holds:
\be\label{eq:CTnD}
\eta(\rho_0,\u_0)=\lamin(\nablaS\u_0)(\x) +\amptwo\rhoav_0(\x) \geq \etac >0, \qquad \etac=\tfrac{1}{2}\amp \etar.
\ee
Then, Euler alignment  \eqref{eq:Int-mono-kinetic},\eqref{eq:V0}--\eqref{eq:CTnD} admits a global smooth solution,  $(\rho(\cdot,t),\u(\cdot,t))\in L^1_+\times W^{1,\infty}(\RR^n)$, and  the following uniform bound holds
\[
|\nabla \u(t,\cdot)|_{L^\infty} \leq \max\left\{|\nabla \u_0|_{L^\infty},\etar,C_0\right\}.
\]
The constant $C_0=C\left(\max_\x\|\nabla \u_0\|, m_0,\phi_+\right)$ is specified  in \eqref{eq:finalM} below.
 \end{theorem}
 
 \noindent
  In the canonical case of metric kernels, Theorem \ref{thm:CS} holds with $\beta_0=0$.  This is summarized in Theorem \ref{thm:CS-metric}.

 \begin{remark}[{\bf Large class of sub-critical initial data}]
 Observe that the threshold condition  \eqref{eq:CTnD} is the same threshold condition which secured decay towards mono-kinetic closure for general class of entropic pressures in  Theorem \ref{thm:threshold}.
 It allows a ``large'' set of sub-critical initial configurations $(\rho_0,\u_0)$   in  \eqref{eq:CTnD} such that   $\lamin(\nablaS\u_0)$ admits negative values 
 $\lamin(\nablaS\u_0)> \etac - \amp \rhoav$ with $ \etac - \amp \rhoav\leq  -\hf\amp\etar$; consult Remark \ref{rmk:what}. One can state  the existence result with a smaller  threshold $\etac$ (thus allowing even larger range of ``negative slopes'' $\lamin(\nablaS\u_0)$) at the expense of more restricted range of initial amplitudes $(\alpha_0,\beta_0)$.
 \end{remark}
 
\noindent
\emph{Proof of Theorem \ref{thm:CS}.} 
Our purpose is to show that the derivatives $\{\pa_ju_i\}$ are uniformly bounded. 
We proceed in four steps.\newline
{\bf Step 1}.  We begin by identifying an invariant region associated with the threshold  $\eta(t,\x)=\lamin(\nablaS\u)+\amptwo\rhoav$. The goal is 
to show that the threshold condition \eqref{eq:CTnD} secures a sub-critical region in  configuration space which persists in time, 
$\eta(\rho_0,\u_0)\geq \etac>0 \ \leadsto \ \eta(\rho(t,\cdot),\u(t,\cdot))\geq \etac>0$.

\noindent
Let ${}^\prime$ abbreviate differentiation along particle path,  
\[
\Box^\prime:=(\partial_t+\u\cdot\nabla)\Box.
\]
Differentiation of \eqref{eq:mono-momentum} implies that  the $n\times n$ velocity gradient matrix, $M=\nabla \u$, satisfies
\bea\label{eq:CSM}
M^\prime + M^2  +\amp \rhoav M =\amp R,
\eea
where  $R$  is the $n\times n$  matrix 
\[
R_{ij}:=\int \frac{\partial \phi}{\partial x_j}(\x,\y)\big(u_i(t,\y)-u_i(t,\x)\big)\rho(t,\y)\dy.
\]
The following observation of the residual $R$ is at the heart of matter;  
in the special case of \emph{metric} kernels it goes back to \cite{CCTT2016} in the 1D case, and  to \cite{HT2017} in the 2D case,
\begin{equation}\label{eq:Rident}
\trace{R} = -(\rhoav)^\prime +  \rhoavpsi, \qquad \psi(t,\x,\y):=\sum_i\left(\frac{\partial \phi}{\partial {x_i}}+\frac{\partial \phi}{\partial {y_i}}\right)u_i(t,\y).
\end{equation}
\emph{Verification of} \eqref{eq:Rident}: integration by parts followed by the mass equation \eqref{eq:EA}${}_1$ yield
\[
\begin{split}
\trace&{R}= \int \sum_i\frac{\partial \phi}{\partial x_i}(\x,\y)\big(u_i(t,\y)-u_i(t,\x)\big)\rho(t,\y)\dy \\
& =-\int \sum_i\frac{\partial \phi}{\partial y_i}(\x,\y)u_i(t,\y)\rho(t,\y)\dy 
- \int \sum_i\frac{\partial \phi}{\partial x_i}(\x,\y)u_i(t,\x)\rho(t,\y)\dy
\\
 & \qquad +\int \sum_i\left(\frac{\partial \phi}{\partial x_i}+\frac{\partial \phi}{\partial y_i}\right)(\x,\y) u_i(t,\y)\rho(t,\y)\dy\\
& = \int  \phi(\x,\y)\nabla_{\y}\cdot(\rho \u)(t,\y)\dy  -\u\cdot\nabla_\x \int \phi(\x,\y)\rho(t,\y)\dy + \int \psi(\x,\y) \rho(t,\y)\dy \\
& = \int  \phi(\x,\y)\times -\rho_t (t,\y)\dy- \u\cdot\nabla_\x \rhoav +  \psi*\rho\\
& = -(\partial_t+\u\cdot\nabla_\x )\rhoav +\rhoavpsi.
\end{split}
\]
We decompose $M$ into its symmetric and skew-symmetric arts, $M=S+\Omega$ and trace the  symmetric part   of \eqref{eq:CSM}, $S=\nablaS\u$, to find
\begin{equation}\label{eq:MS}
S^\prime + S^2 +\Omega^2  + \amp\rhoav S  = \amp R_{S}, \qquad R_{S}:=\hf(R+R^\top).
\end{equation}
Since $\Omega$ is skew-symmetric then $S^\prime +S(S+\amp\rhoav{\Bbb I})\geq \amp R_S$.
Hence, the minimal eigen-pair of $S$, $\lamin:=\min_\lam \lam(S)$ with the corresponding unit eigenvector $\w_-$,  satisfies
\be\label{eq:X}
\lamin^\prime(t,\x) +\lamin(t,\x)\big(\lamin(t,\x)+\amp\rhoav(t,\x)\big) \geq \amp \langle R_{S}\w_-,\w_-\rangle.
\ee
Here comes the main motivation  for \eqref{eq:Rident}: we use it to express \eqref{eq:X}  in terms of $\eta:=\lambda_-+\amp\rhoav$ as follows:
\[
\eta'+ (\eta-\amp\rhoav)\eta \geq \amp \chi, \qquad \chi:=\langle (R_{S}-\trace{R_S})\w_-,\w_-\rangle +\psi*\rho,
\]
and since by uniform thickness $\rhoav\geq \etar$ hence,
\begin{equation}\label{eq:lowereta}
\eta^\prime(t,\x) \geq \eta(\amp\etar-\eta)-\amp|\chi|. 
\end{equation}
Our purpose is to secure that the quantity on the right of \eqref{eq:lowereta} is positive for the threshold $\eta\downarrow \etac$, so that $\eta(t,\cdot)$ always increases whenever it approaches $\etac=\tfrac{1}{2}\amp\etar$ from above, and hence $\eta(\rho(t,\cdot),\u(t,\cdot))\geq \etac$ remains an invariant region in time. To this end, we need to upper-bound $|\chi|$.

The rank-one integrand in $R$ implies 
$\displaystyle |\langle R_S\w_-,\w_-\rangle|\leq |\nabla_\x\phi|_\infty \delta u(t)m_0$,
 and the velocity fluctuations bound \eqref{eq:noincrease}${}_2$ yields\footnote{The entries integrated in $R_{ij}\mapsto r_is_j$ with $r_i=\phi_{x_i}$ and $s_j=u_i(t,\y)-u_i(t,\x)$ yields,
$\langle {\mathbf r},\bw_-\rangle \langle{\mathbf s},\bw_-\rangle   -\sum r_is_i \leq   2 |\nabla_{\x}\phi|\delta u(t)$.}
\begin{equation}\label{eq:xyz}
\big|\big\langle (R_S-\trace{R})\w_-,\w_-\big\rangle\big| \leq  2\alpha_0 m_0
\qquad \alpha_0=|\nabla_\x\phi|_\infty \delta u_0;
\end{equation}
the kinetic energy bound \eqref{eq:noincrease}${}_3$ yields
\[
|\rhoavpsi|\leq \beta_0 m_0, \qquad \beta_0= |(\nabla_\x+\nabla_{\y})\phi|_\infty \sqrt{\frac{2\kE_0}{m_0}}. 
\]
Therefore, since the initial amplitudes assumed in \eqref{eq:V0} are not too large,
\[
\amp|\chi|\leq \amp|\langle (R_{S}-\trace{R_S})\w_-,\w_-\rangle|+\amp|\psi*\rho| \leq \tfrac{1}{4}(\amp\etar)^2,
\]
and \eqref{eq:lowereta} yields
\[
\eta^\prime(t,\x) >  \eta(\amp\etar-\eta)-\tfrac{1}{4}(\amp\etar)^2_{|\eta=\etac}=0, \qquad \etac=\tfrac{1}{2}\amp\etar. 
\]
Thus,  the initial threshold bound $\eta(\rho_0,\u_0)\geq \etac$ with $\etac=\tfrac{1}{2}\amp\etar$, will persist for all time,  $\eta(\rho,\u)\geq \etac$.

\noindent
{\bf Step 2}. Next we upper-bound the skew-symmetric part of \eqref{eq:CSM}
\[
\Omega' + \hf(S\Omega+\Omega S) +\amp\rhoav \Omega = \amp R_\Omega, \qquad R_\Omega:=\tfrac{1}{2}\left(R-R^\top\right).
\]
If $(\mumax,\z_+)$ is the  maximal eigen-pair of $\Omega$ with purely imaginary eigenvalue  $\mumax$ such that $ |\mumax|=\max_\mu |\mu(\Omega)|$ and normalized eigenvector $\z_+$, then
\[
\mumax'(t,\x)+\mumax(t,\x)\big(\langle S\z_+,\z_+\rangle + \amp\rhoav\big)
=\amp \langle R_\Omega\z_+,\z_+\rangle.
\]
The threshold established in step 1  tells us,  $\langle S\z_+,\z_+\rangle + \amp\rhoav \geq \eta_c >0$ and hence
\[
|\mumax|'(t,\x) + \eta_c|\mumax(t,\x)|\leq |\mumax|'(t,\x)+|\mumax(t,\x)|\big(\langle S\z_+,\z_+\rangle + \amp\rhoav\big) \leq 
\amp |\langle R_\Omega\z_+,\z_+\rangle|.
\]
As before,  $|\langle R_\Omega\z_+,\z_+\rangle|\leq \alpha_0m_0< \tfrac{1}{8}\amp\etar^2$ and we end up with a uniform bound of the vorticity whch involves the constant ${\tfrac{1}{8}\amp^2\etar^2}/{\etac}=\tfrac{1}{4}\amp\etar$,
\[
 \|\Omega(t,\x)\| \leq \gamma_0, \qquad \gamma_0:=\max\left\{\max_\x \|\Omega_0(\x)\|, \tfrac{1}{4}\amp\etar\right\}.
\]

\noindent
{\bf Step 3}. Now we can bound $\|S(t,\x)\|$, that is, the maximal eigenvalue of $\lamax=\max_\lam \lam(S)$. We revisit \eqref{eq:MS}, this time with the upper bound of $-\Omega^2$, to derive the reverse inequality \eqref{eq:lowereta} for $\zeta(t,\x):=\lamax(t,\x)+\amp\rhoav(t,\x)$, 
\[
\zeta^\prime(t,\x) +\amp\etar \zeta(t,\x) \leq \zeta^\prime +\amp\rhoav \zeta\leq \amp|\langle R_S\bw_+,\bw_+\rangle|+\|\Omega\|^2.
\]
Again, since  the initial amplitudes assumed in \eqref{eq:V0} are not too large, $\amp|\langle R_S\bw_+,\bw_+\rangle|\leq \tfrac{1}{8}(\amp\etar)^2$  and together with $\|\Omega\|\leq \gamma_0$, this  yields
\[
\max_\x \lamax(t,\x) \leq \max\Big\{\max_\x\lamax(S_0)(\x), \delta_0\Big\}, \qquad \delta_0=\tfrac{1}{8}\amp\etar + \frac{\gamma_0^2}{\amp\etar}..
\]
And finally, combined with the lower threshold  $\lamin(S)\geq -\amp\rhoav \geq -\amp\phi_+m_0$ we conclude 
\[
\|S(t,\x)\|\leq \max\left\{\max_\x \lamax(S_0)(\x),  \delta_0, \amp\phi_+m_0\right\}.
\] 
\noindent
{\bf Step 4}. The bounds of $S$ and $\Omega$ imply that $\nabla\u$ is uniformly bounded in time 
\begin{equation}\label{eq:finalM}
\left|\frac{\partial u_i}{\partial x_j}(t,\x)\right| \leq \max\left\{ \max_\x\|\nabla \u_0\|, \frac{1}{4}\amp\etar, \delta_0,
\amp\phi_+m_0\right\}. \qquad \square
\end{equation}

\begin{remark}
Observe that the velocity fluctuations $\delta u(t)$ decay exponentially in 
time, see \eqref{eq:exp-fluc-vel-in-time} below. If we use this improved bound in \eqref{eq:xyz} then one can deduce an improved threshold  condition with larger set of restricted initial fluctuations, e.g., \cite[Remark 6.1]{ST2020a}.
\end{remark}
\medskip

\subsection{Uniform thickness and flocking}\label{sec:thickness} 
We now turn to the question of uniform thickness assumed in Theorem \ref{thm:CS} which  is addressed in the next sections for the two main classes of heavy-tailed and short-range kernels.

\vspace*{-0.3cm}
\subsubsection{Uniform thickness with heavy-tailed kernels}\label{sec:heavy-thickness} 
Consider communication kernels, quantified in terms of the Pareto-type tail
\begin{equation}\label{eq:Pareto}
\phi(\x,\y)\geq C \langle|\x-\y|\rangle^{-\theta}, \qquad \theta\in (0,1)
\end{equation}
A key feature of such heavy-tailed kernels is that their alignment dynamics maintains global communication:  each  part of the crowd  with mass distribution $\rho(\x)\dx$ communicates \emph{directly} with every other part with mass distribution $\rho(\y)\dy$. 
Indeed, in this case,
\[
\phi_-(t):=\min_{\x,\y\in\suppr} \phi(\x,\y)\gtrsim \langle D(t)\rangle^{-\theta}>0, \qquad D(t):=\textnormal{diam}(\suppr)
\]
and for $\theta<1$ this implies that $\phi_-(t)$ remains uniformly bounded in time away from zero. 
The uniform-in-time bound follows by combining two standard arguments:
 
\#1. Decay of velocity fluctuations. In the case of mono-kinetic closure, $\presure\equiv 0$, the momentum equation \eqref{eq:EA}$\mbox{}_2$ decouples into $n$ scalar equations, each of which satisfies a maximum principle; in fact there is a decay quantified in  \cite[Theorem 2.1]{TT2014},\cite[Theorem 1.1]{HT2017},
\begin{equation}\label{eq:max-principle}
\ddt \delta u(t) \leq -(\amp m_0\phi_-(t))\delta u(t), \qquad 
\delta u(t) =\max_{\bx,\by\in \suppr} |\bu(t,\bx)-\bu(t,\by)|.
\end{equation}

\#2. The decay rate of $\delta u(t)$,
\begin{subequations}\label{eqs:D-and-delta}
\begin{equation}\label{eq:delta-by-D}
\ddt \delta u(t) \leq -C\amp m_0 \langle D(t)\rangle^{-\theta}\delta u(t)
\end{equation}
is dictated by the dispersion of $\suppr$, which in turn does not exceed
\begin{equation}\label{eq:D-by-delta}
\ddt D(t) \leq \delta u(t).
\end{equation}
\end{subequations}
It follows \cite{HL2009} that $H(t):=C\amp m_0\langle D(t)\rangle^{1-\theta}+(1-\theta)\delta u(t)$ is non-increasing, $\dot{H}(t) \leq 0$,
and therefore $\suppr$ is kept uniformly bounded in time, 
\begin{equation}\label{eq:finite-dispersion}
 D(t)\leq D_\infty:=\left(\frac{H_0}{\amp m_0}\right)^{\tfrac{1}{1-\theta}}, \qquad H_0=C\amp m_0\langle D_0\rangle^{1-\theta}+(1-\theta)\delta u_0.
\end{equation}
In particular, $\phi_-(t) \geq C \langle D(t)\rangle^{-\theta} \geq C \langle D_\infty\rangle^{-\theta}>0$ and  uniform thickness \eqref{eq:Int-thickness} follows
\begin{equation}\label{eq:heavy-thickness}
\rhoav(t)\geq \phi_-(t)\int \rho(t,\x)\dx \geq \etar:=C\langle D_\infty\rangle^{-\theta}m_0.
\end{equation}
An alternative derivation of the uniform-in-time bounds is outlined in \S\ref{sec:uniform} below.

\begin{remark}[{\bf Flocking I}]
The dispersion bound \eqref{eq:finite-dispersion} implies an exponential decay of velocity fluctuations  
\begin{equation}\label{eq:exp-fluc-vel-in-time}
\delta u(t) \lesssim e^{-\amp m_0 D_\infty^{-\theta}t}\delta u_0.
\end{equation}
Since the   mean velocity is time-invariant, $\displaystyle \overline{\u}(t):=\frac{1}{m_0}\int \rho\u(t,\x)\dx$, there follows the flocking behavior
\begin{equation}\label{eq:utoubar} 
\mathop{\max}_{\x\in \suppr}|\u(t,\x)-\overline{\u}_0|_\infty \lesssim e^{-\ampzero D^{-\theta}_\infty t} |\delta \u_0|_\infty. 
\end{equation}
Flocking behavior for heavy-tailed metric-based kernels $\phi(\x,\y)\mapsto \phi(|\x-\y|)$  goes back to Cucker-Smale \cite{CS2007,HT2008,HL2009}. Here we observe that it extends to general heavy-tailed symmetric kernels. 
It corresponds to the flocking behavior at the level of agent-based description e.g., \cite[definition 1.1]{MT2011}, in which a cohesive flock of a finite diameter $\max_{i,j}|x_i(t)-x_j(t)|\leq \diam_\infty<\infty$, is approaching a limiting velocity, $\max_{i,j}|\v_i(t)-\v_j(t)|\rightarrow 0$ as $t\rightarrow \infty$.
\end{remark}
\begin{remark}[{\bf Flocking II}]
The  enstrophy in \eqref{eq:Int-Efluc}  drives the  the decay of the energy fluctuations 
\[
\begin{split}
\ddt \delE(t) \leq & -\frac{C\ampzero}{\langle D(t)\rangle^{\theta}}\delE(t), \\
 & \quad \delE(t):=
\left(\iint \Big(\hf|\u(\x)-\u(\y)|^2 +e(\x)+e(\y)\Big)\d\rho(\x)\d\rho(\y)\right)^{\nicefrac{1}{2}}.
\end{split}
\]
Observe that the last inequality  applies to general  Euler alignment systems \eqref{eq:EA}, \emph{independent} of the specific closure for the pressure tensor. In particular, \underline{if} we can control the dispersion
\begin{equation}\label{eq:D-dispersion}
D(t)\lesssim \langle t\rangle^{\gamma}, \quad \gamma\geq 0,
\end{equation}
then we would conclude that both ---  the $L^2$-velocity fluctuations  and internal energy fluctuations decay to zero
\begin{equation}\label{eq:expalign}
\delE(t) \lesssim  e^{-\eta t^{1-\gamma\theta}}\delE(0),
\qquad \eta=\frac{C\ampzero}{2(1-\gamma\theta)u_\infty}, \quad \gamma\theta<1.
\end{equation}
We mention three examples.\newline
1 (A uniformly bounded velocity). $|\u(t,\bx)|\leq u_\infty$ implies --- in view of \eqref{eq:D-by-delta}, that \eqref{eq:D-dispersion} holds with $\gamma=1$, $D(t)\leq D_0+2u_\infty t$, and $L^2$-flocking follows for $\theta<1$.
This covers  the mono-kinetic scenario discussed in the previous remark.
\newline
2 (Matrix communication kernels). When \eqref{eq:EA} is driven by symmetric positive definite \emph{matrix} kernel, $\phi\mapsto \Phi(\x,\y)$, then a dispersion bound is secured with $\gamma=\nicefrac{2}{3}$, and \eqref{eq:expalign}  follows for $\theta<2/3$, \cite[Proposition 3.1]{ST2021},\cite[Appendix C]{Tad2023}.\newline
3 (Singular kernels). When $\phi(\x,\y)\sim |\x-\y|^{-\alpha}, \ \alpha\in(n,n+2)$, their enstrpohy enforces that \eqref{eq:D-dispersion} holds with $\gamma=\gamma_{\alpha,n}>1$ \cite[Appendix E]{Tad2023}, and flocking follows for a restricted set of $\theta<1/\gamma$. This proves flocking  independent of the specific closure for the pressure and independent whether the kernel $\phi(\x,\y)$ is long- or short-range.
\end{remark}

\vspace*{-0.3cm} 
\subsubsection{Uniform thickness with short-range kernels}\label{sec:short-thickness}

We restrict attention to compactly supported metric kernels \eqref{eq:metricphi}.
In this case,  alignment takes place in local neighborhoods of size $\leq D_0$, which is assumed much smaller than the diameter of the ambient space --- the $2\pi$-periodic torus $\Omega={\Bbb T}^n$. 
With lack of global, direct communication, existence and large-time behavior of strong solution of \eqref{eq:EA} depends on having a bounded, non-vacuous density, 
\begin{equation}\label{eq:non-vacuous}
0<\rho_- \leq \rho(t,\x) \leq \rho_+ <\infty. 
\end{equation}
In particular, uniform thickness \eqref{eq:Int-thickness} follows
\begin{equation}\label{eq:short-thickness}
\rhoav(t)\geq \etar:= \int \phi(|\x|)\dx \, \rho_- >0.
\end{equation}\begin{remark}[{\bf Flocking III}]
We mention the flocking of alignment dynamics with short-range kernels in finite tours, as long as the the non-vacuous condition \eqref{eq:non-vacuous} holds, \cite[Theorem 3]{Tad2021},\cite[\S4.4]{Shv2024}. Again, flocking of non-vacuous dynamics holds independent of the closure of pressure: the key question is securing  the uniform-in-time bound \eqref{eq:short-thickness} or, as noted in \cite[Remark p. 499]{Tad2021}  at least $\displaystyle \int^\infty \rho_-(t)\dt =\infty$ which in tun would lead to (non-uniform) thickness, 
$\displaystyle \int^\infty \min_\x \rhoav(t,\x)\dt=\infty$.
This question was addressed in the case of 1D singular kernel in \cite{ST2017a,DKRT2018}, but is open for bounded short-kernels.
\end{remark}

\subsection{Uniform bounds revisited}\label{sec:uniform} A  necessary main ingredient in the analysis of \eqref{eq:EA} is the uniform-in-time bounds of $\text{diam}(\suppr), \phi_-(t)$, and the amplitude of velocity $\displaystyle \mathop{\max}_{\x\in \suppr}|\u(\bx,t)|$. An alternative approach to the standard arguments in \S\ref{sec:heavy-thickness} is advocated in our work \cite[Lemma 3.2]{ST2020a}. Here, we extend our argument to general symmetric kernels. The next lemma shows that whenever one has a uniform bound of $|\u(\bx,t)|+|\x|$ for the \emph{restricted} class of lower-bounded $\phi$'s which scales like ${\mathcal O}(1/\min \phi)$, then it implies a uniform bound of 
$|\u(\bx,t)|+|\x|$ for the general class of admissible $\phi$'s \eqref{eq:Pareto}. 

\begin{lemma}[{\bf The reduction to lower-bounded $\phi$'s}]\label{cor_reduce}
Consider \eqref{eq:EA} with a with the restricted class of uniformly lower-bounded $\phi$'s:
\begin{equation}\label{phiminus}
0< \phi_- \le \phi(\bx,\by) \le \phi_+ < \infty.
\end{equation}
Assume that the solutions  $(\widetilde{\rho},\widetilde{\bu})$ associated with the restricted \eqref{eq:EA},\eqref{phiminus},  satisfy the uniform bound (with constants   $C_\pm$ depending on $\phi_+,m_0$ and $\kE_0$)
\begin{equation}\label{eq:ubound}
\max_{t\ge 0,\, \x\in\supp\widetilde{\rho}(\cdot,t)} (|\widetilde{\u}(\x,t)|+|\x|) \le \max\left\{\Cplus\cdot\hspace*{-0.4cm}\max_{\quad \x\in\textnormal{supp }\widetilde{\rho}_0} (|\widetilde{\u}_0(\x)| + |\bx|)\,,\frac{\Cminus}{\phi_-}\right\}.
\end{equation}
Then the following holds for solutions associated with a general class of  admissible kernels $\phi$ satisfying \eqref{eq:Pareto}: if $(\rho,\bu)$ is a smooth solution of \eqref{eq:EA}, then there exists $\beta>0$ depending on the initial data  $(\rho_0,\bu_0)$,  such that $(\rho,\bu)$ coincides with  the solution, $(\widetilde{\rho}_\beta,\widetilde{\bu}_beta)$, associated with the lower-bounded  
\[
\phi_\beta(\bx,\by) := \max\{\phi(\bx,\by),\beta\}
\]
\end{lemma}
This means that if $\phi$ belongs to the general class of admissible kernels \eqref{eq:Pareto}, then we can  assume, without loss of generality, that $\phi$ coincides with the lower bounded $\phi_\beta$ and hence the uniform bound \eqref{eq:ubound} holds with $\phi_-=\beta$. The justification of this reduction step is outlined below.

\medskip\noindent
\emph{Proof of Lemma \ref{cor_reduce}.}
By  \eqref{eq:Pareto} one could take large enough $r$ such that
$\displaystyle  \mathop{r \cdot \min}_{\ \ \ |\x-\y|\leq r}\phi(\x,\y) \ge 2\Cminus$ and 
\begin{equation}\label{r0large}
r\ge 2\Cplus\cdot \hspace*{-0.4cm}\max_{\quad \bx\in\textnormal{supp }\rho_0}\hspace*{-0.4cm}(|\bu_0(\bx)|+|\bx|).
\end{equation}
Let $\displaystyle \beta:=\min_{|\x-\y|\leq r}\phi(\x,\y)$. By assumption, (\ref{eq:ubound}) holds for the lower-bounded $\phi_\beta$, so that
\begin{equation}
\max_{t\ge 0,\, \bx\in\textnormal{supp }\rho_\beta(\cdot,t)} (|\bu_\beta(\bx,t)|+|\bx|) \le \max\Big\{\Cplus\cdot\hspace*{-0.4cm}\max_{\quad \bx\in\textnormal{supp }\rho_0}\hspace*{-0.4cm} (|\bu_0(\bx)|+|\bx|), \frac{\Cminus}{\beta}\Big\}
\end{equation}
where $(\rho_\beta,\bu_\beta)$ is the smooth solution of \eqref{eq:EA} with interaction kernel  $\phi_\beta$, which we assume to exist. 
By definition,
\begin{equation}\label{eq:rtwo}
\frac{\Cminus}{\beta} = \frac{\Cminus}{\displaystyle \min_{|\x-\y|\leq r}\phi(\x,\y)} \le \frac{ r}{2}.
\end{equation}
Fix $t\ge 0$, then (\ref{r0large})--\eqref{eq:rtwo} imply that the distance for any $\bx,\by \in \textnormal{supp }\rho_\beta(\cdot,t)$ does not exceed 
\begin{equation}
|\bx-\by| \le |\bx|+|\by| \le 2\max\Big\{\Cplus\cdot \hspace*{-0.4cm}\max_{\quad\bx\in\textnormal{supp }\rho_0}\hspace*{-0.4cm} (|\bu_0(\bx)|+|\bx|), \frac{\Cminus}{\beta}\Big\} \leq r.
\end{equation}
Thus, for any  $\bx,\by \in \textnormal{supp }\rho_\beta(\cdot,t)$ there holds $\phi(\bx,\by)\geq   \beta$ and hence  $\phi_\beta$ coincides with $\phi$ for $(\bx,\by) \in \suppr$,
so that the dynamics of $(\rho_\beta,\bu_\beta)$ coincides with $(\rho,\bu)$. \hfill $\square$

\noindent
Example. If the uniform lower bound $\phi(\bx,\by)\geq \phi_-$ holds then 
according to \eqref{eq:max-principle}
\[
\delta u(t)\leq \delta u_0 \cdot e^{-(\ampzero\phi_-)t},
\]
and hence $\displaystyle D(t)\leq D_0 + \frac{1}{\ampzero\phi_-}$.
Therefore, existence of strong solutions and their flocking behavior follows  long range $\phi$'s with Pareto's tail \eqref{eq:Pareto}.


\begin{thebibliography}{1}

\bibitem[Ave2003]{Ave2003}
J. Avery, Information Theory and Evolution, World Scientific, 2003.
 
\bibitem[CCR2009]{CCR2009}
J. A. Canizo, J.A. Carrillo and J. Rosado, 
A well-posedness theory in measures for 
kinetic models of collective motion, 
Math. Mod. Meth. Appl. Sci., 21 (2009), 515-539.

\bibitem[CFTV2010]{CFTV2010}
    J.~A. Carrillo, M.~Fornasier, G.~Toscani and F.~Vecil.
    \newblock Particle, kinetic, and hydrodynamic models of swarming, 2010.
    \newblock in ``\emph{Mathematical Modeling of Collective Behavior 
    in Socio-Economic and Life Sciences}'' 
    (G. Naldi L. Pareschi G.  Toscani eds.), Birkhauser, pages 297--336, 2010.


\bibitem[CCP2017]{CCP2017}
J. Carrillo, Y.-P. Choi and S. Perez,
A review on attractive-repulsive hydrodynamics
for consensus in collective behavior,
in ``\emph{Active Particles,  Volume 1. 
Advances in Theory, Models, and Applications}''
(N. Bellomo, P. Degond and E. Tadmor, eds.),
Birkhäuser 2017.


\bibitem[CCTT2016]{CCTT2016}
J. A. Carrillo, Y.-P. Choi, E. Tadmor and C. Tan,
Critical thresholds in 1D Euler equations with nonlocal forces,
Mathematical Models and Methods in Applied Sciences 26(1) (2016) 185-206.

\bibitem[CFRT2010]{CFRT2010}
J. Carrillo, M. Fornasier, J. Rosado and G. Toscani, 
Asymptotic flocking dynamics for the kinetic Cucker-Smale model, 
SIAM J. Math. Anal. 42(218) (2010), 218-236.

\bibitem[Cha1978]{Cha1978}
G. Chaitin, {Towards a mathematical definition of ``Life''}, in ``\href{http://home.thep.lu.se/~henrik/mnxa09/Chaitin1979.pdf}{The
Maximum Entropy Formalism}'' (R. D. Levine and M. Tribus, eds.), MIT Press,
1979, pp. 477--498

\bibitem[CS2007]{CS2007}
{F. Cucker and S. Smale}, 
{Emergent behavior in flocks},
IEEE Trans. Autom. Control, 52, no. 5, (2007): 852-862.

\bibitem[Daf2005]{Daf2005}
C. Dafermos, Hyperbolic conservation laws in continuum physics. Vol. 3. Berlin, Springer, 2005.


\bibitem[DKRT2018]{DKRT2018}
T. Do, A. Kiselev, L. Ryzhik and C. Tan, Global regularity for the fractional Euler alignment system.
Arch. Ration. Mech. Anal. 228(1), 1–37 (2018)

\bibitem[ELT2001]{ELT2001}
S. Engelberg, H. Liu and E. Tadmor, 
Critical thresholds in Euler-Poisson equations,
Indiana Univ. Math J.  50 (2001), 109-157


\bibitem[FP2024]{FP2024}
M. Fabisiak and J. Peszek,  Inevitable monokineticity of strongly singular alignment. \href{https://doi.org/10.1007/s00208-023-02776-7}{Math. Ann. 390, (2024) 589-637.}

\bibitem[FK2019]{FK2019}
A. Figalli and M.-J. Kang,
\newblock A rigorous derivation from the kinetic Cucker–Smale model to the pressureless Euler system with nonlocal alignment
\newblock Anal. PDE, 1293)  (2019), 843-866.

\ifx
\bibitem[Gia2008]{Gia2008}
I. Giardina,
\newblock Collective behavior in animal groups: Theoretical models and empirical studies
\newblock HFSP Journa 2(4) (2008) 205-219.
 \fi
 
  \bibitem[God1962]{God1962}
 S. K. Godunov, The problem of a generalized solution in the theory of quasilinear equations and in gas dynamics,  Russian Mathl Surveys 17(3) (1962): 145.
 
\bibitem[HL2009]{HL2009}
S.~Y. Ha and J.~G Liu,
\newblock A simple proof of the {Cucker-Smale} flocking dynamics and mean-field
  limit.
\newblock {Comm. in Math. Sciences}, 7(2):297--325, 2009.

\bibitem[HT2008]{HT2008}
{S.-Y. Ha and E. Tadmor},
{From particle to kinetic and hydrodynamic descriptions of flocking},
Kinetic and Related Models, 1(3), (2008), 415-435.

\bibitem[HT2017]{HT2017}
S. He and E. Tadmor,
\newblock Global regularity of two-dimensional flocking hydrodynamics.
\newblock Comptes rendus - Mathématique Ser. I 355 (2017) 795-805.

\ifx
\bibitem[HT2019]{HT2019}
S. He and E. Tadmor,
\newblock A game of alignment: collective behavior of multi-species.
\newblock ArXiv:1908.11019 (2019)
\fi

\bibitem[Kru1970]{Kru1970}
S. Kruzkov, First order quasilinear equations in several independent variables, Math. of the USSR-Sbornik 10.2 (1970): 217.

\bibitem[Lax1957]{Lax1957}
P. D. Lax, Hyperbolic systems of conservation laws. Comm. Pure Appl. Math. 10 (1957), 537-566.

\bibitem[Lax1971]{Lax1971}
P.  Lax, Shock waves and entropy, in ``\emph{Contributions to nonlinear functional analysis}'', Academic Press, 1971. 603-634.

\ifx
\bibitem[Lax1973]{Lax1973}
P. Lax,  Hyperbolic systems of conservation laws and the mathematical theory of shock waves, {SIAM}, 1973.
\fi

\ifx
\bibitem[LL2013]{LL2013}
{Y. Lee and H. Liu}, 
{Thresholds in three-dimensional restricted Euler-Poisson equations},
Physica D,  262, (2013), 59-70.
\fi

\bibitem[LS2020]{LS2020}
D. Lear and R. Shvydkoy,
\newblock
 Existence and stability of unidirectional flocks in hydrodynamic Euler alignment
systems,
\newblock Analysis \& PDE 15(1) (2022), 175-196.


\bibitem[Les2020]{Les2020} 
T. M. Leslie, 
\newblock On the Lagrangian Trajectories for the One-Dimensional Euler Alignment Model without Vacuum Velocity,
\newblock  Comptes Rendus Mathématique 358(4) (2020) 421-433.

\bibitem[LLST2022]{LLST2022}
D. Lear, T. M. Leslie, R. Shvydkoy and E. Tadmor, 
Geometric structure of mass concentration sets for pressureless Euler alignment systems, Advances in Mathmematics 401(4) (2022) 108290.

\ifx
\bibitem[LT2001]{LT2001}
H. Liu and E. Tadmor,
Critical thresholds in a convolution model for nonlinear conservation laws,
SIAM Journal on Mathematical Analysis 33 (2001), 930-945
\fi

\bibitem[LT2002]{LT2002}
H. Liu and E. Tadmor,
Spectral dynamics of the velocity gradient field in restricted flows,
Communications in Mathematical Physics 228 (2002), 435-466.

\ifx
\bibitem[LT2003]{LT2003}
H. Liu and E. Tadmor, 
Critical thresholds in 2D restricted Euler-Poisson equations,
SIAM Journal of Applied Mathematics 63 (2003) 1889-1910.


\bibitem[LT2004]{LT2004}
H. Liu and E. Tadmor, 
Rotation prevents finite-time breakdown,
Physica D 188 (2004) 262-276.


\bibitem[MOA2010]{MOA2010}
N.Mecholsky, E. Ott and T. Antonsen,
Obstacle and predator avoidance in a model for flocking,
Physica D, 239(12) 2010, 988-996.
\fi

\bibitem[MMPZ2019]{MMPZ2019}
P. Minakowski, P.B.  Mucha, J. Peszek and E. Zatorska,  Singular Cucker-Smale dynamics. in ``\emph{Active Particles, Volume 2: Advances in Theory, Models, and Applications}'' (N. Bellomo, P. Degond, E. Tadmor eds.),  Springer 2019, pp. 201-243.

\bibitem[MT2011]{MT2011}
S. Motsch and E. Tadmor,
A new model for self-organized dynamics and its flocking behavior,
Journal of Statistical Physics 144(5) (2011) 923-947.

\bibitem[MT2014]{MT2014}
S. Motsch and E. Tadmor, 
Heterophilious dynamics enhances consensus,
SIAM Review 56(4) (2014) 577-621.

\bibitem[NP2022]{NP2022}
R. Natalini and T. Paul, On the mean field limit for Cucker-Smale models. Discrete
Contin. Dyn. Syst. B. 27 (5) (2022) 2873-2889.


\bibitem[NS2022]{NS2022}
V. Nguyen and R. Shvydkoy,
Propagation of chaos for the Cucker-Smale systems under heavy tail communication, Comm.  in PDEs 47(9) (2022) 1883--1906.

\bibitem[PT2022]{PT2022}
T. Paul and E. Tr\'{e}lat,
From microscopic to macroscopic scale equations: mean field, hydrodynamic and graph limits, arXiv:2209.08832. (2022)

 \bibitem[Sch1944]{Sch1944} 
 E. Schr\" {o}dinger,  What is Life -- the Physical Aspect of the Living Cell. Cambridge
University Press. ISBN 978-0-521-42708-1.
 
\bibitem[ST2020a]{ST2020a} 
R. Shu and E. Tadmor, 
\newblock Flocking hydrodynamics with external potentials,
\newblock  Archive Rational Mech. and Anal. 238 (2020) 347-381.   

\bibitem[Shv2024]{Shv2024} R. Shvydkoy, Environmental Averaging. EMS Surveys in Math Sci. 11 (2024), 277-413

\bibitem[Shv2025]{Shv2025}
R. Shvydkoy, Global well-posedness and relaxation for solutions of the Fokker-Planck-Alignment equations, arXiv:2412.20294v2.

\bibitem[ST2017a]{ST2017a}
R. Shvydkoy and E. Tadmor,
Eulerian dynamics with a commutator forcing,
\href{https://doi.org/10.1093/imatrm/tnx00}{Trans. in Math. and Applications, 1(1) (2017), 1-26}.

\bibitem[ST2017b]{ST2017b}
R. Shvydkoy and E. Tadmor,
Eulerian dynamics with a commutator forcing II: flocking, Discrete and Continuous Dynamical Systems-A 37(11) (2017) 5503-5520.

\bibitem[ST2018]{ST2018}
R. Shvydkoy and E. Tadmor,  Eulerian dynamics with a commutator forcing III. Fractional diffusion of
order $0<\alpha<1$. Phys. D 376/377, 131-137 (2018)


\bibitem[ST2020b]{ST2020b} 
R. Shvydkoy and E. Tadmor,
\newblock Topologically-based fractional diffusion and emergent dynamics with short-range interactions,
\newblock  SIAM J. Math. Anal. 52(6) (2020) 5792-5839.

\bibitem[ST2021]{ST2021} 
R. Shu and E. Tadmor, 
\newblock Anticipation breeds alignment,
\newblock Archive Rational Mech. and Anal.  240 (2021) 203-241.

\bibitem[Tad1986]{Tad1986}
E. Tadmor, 
A minimum entropy principle in the gas dynamics equations
Applied Numerical Mathematics 2 (1986), 211-219.

\bibitem[Tad2021]{Tad2021}
E. Tadmor, On the mathematics of swarming: emergent behavior in alignment dynamics, 
Notices of the AMS 68(4) (2021) 493-503.

\bibitem[Tad2023]{Tad2023} 
E. Tadmor, Swarming: hydrodynamic alignment with pressure, Bulletin AMS 60(3) (2023) 285-325.

\bibitem[TT2014]{TT2014}
  E. Tadmor and C. Tan, Critical thresholds in flocking hydrodynamics with non-local alignment, Philosophical Transactions of the Royal Society A: Math., Phys. and Engin. Sciences 372.2028 (2014): 20130401.

\ifx
\bibitem[Tam2018]{Tam2018}
I. Tamvakis, Ioannis, Quantifying life (https://www.researchgate.net/publication/33687, 2018.
\fi

\bibitem[VZ2012]{VZ2012}
{T.~Vicsek and A. Zefeiris},
{Collective motion}, Physics Reprints, {517} (2012) 71-140.

\end{thebibliography}
\end{document}